\documentclass{article}
\usepackage{amsmath}
\usepackage{amsfonts}
\usepackage{amssymb}
\usepackage{tikz}
\usepackage{xcolor}        % for text color
\usepackage{pagecolor}     % for background color
\usepackage[utf8]{inputenc}
\usepackage{cite}

\usetikzlibrary{graphs}
\newtheorem{lemma}{Lemma}
\newtheorem{Theorem}{Theorem}

\newcommand{\beq}{\begin{equation}}
\newcommand{\eq}{\end{equation}}

\newcommand{\smallSq}[2]{\left[ \begin{smallmatrix} #1 \\ #2 \end{smallmatrix} \right]}

\title{Generalization and Power of Kocay's Lemma in Graph Reconstruction}
\author{Peter Stark}

\begin{document}
\maketitle	
\begin{abstract}
This paper generalizes Kocay's lemma, with particular applications to graph reconstruction, as well as discussing and proving aspects around the power of these generalizations and Kocay's original lemma, with a result on the reconstruction of the multiplicity of a tree $T$ as a subgraph of $G$.
\end{abstract}
\section{Introduction}

A graph, $G$, is a finite set of vertices $V(G)$, with a set $E(G)$ of pairs of those vertices. We say $G \cong H$ when there is bijection $f: V(G) \rightarrow V(H)$ such that, for any $v_1, v_2 \in V(G)$, $(v_1,v_2) \in E(G)$ iff $f(v_1), f(v_2) \in E(H)$. 

We then say that $G$ and $H$ are hypomorphic when there's a bijection $r: V(G) \rightarrow V(H)$ such that $G-v \cong H-r(v)$. The Reconstruction Conjecture, or the Kelly-Ulam Conjecture \cite{kelly1957} is the statement that hypomorphism is equivalent to isomorphism, or that two graphs are isomorphic iff they have the same deck. A property is called reconstructible or recognizable when it can be deduced from the deck of a graph.

An equivalent statement defines the deck as the multiset of isomorphism classes of single-vertex removed subgraphs.

One of the most powerful tools for reconstruction (which also contains all the information of the deck) is Kelly's lemma, which states that the count of any subgraph with fewer than $n$ vertices, induced or not, are reconstructible. 
 
This is already an extremely useful tool for graph reconstruction. Not only does it allow the reconstruction of the deck, making it equivalent to knowing the deck, but it also gives explicit and deep information about many smaller factors of graphs, leaving it as an essential tool in graph reconstruction. 

Beyond this, though, there's a further generalization. Kocay's lemma allows the explicit reconstruction of sums of counts of graphs with n vertices, so long as they are all unions of graphs with fewer than n vertices. This paper extends and generalizes Kocay's lemma.

\section{Notation}

%Multisets are denoted by square brackets, and $[A*m: P(A)]$ refers to the multiset of all $A$ satisfying $P(A)$ with multiplicity of $m$ satisfying $P(A)$. $Mult(x \in S)$ refers to how many times $x$ appears in $S$. Given a multiset $S$, let $set(S)$ be the set of elements in $S$.

%$\smallSq{G}{H}$ refers to the number of induced subgraphs isomorphic to $H$ in $G$. If we let $ind(P)$ be 1 when $P$, 0 otherwise, then $\smallSq{G}{H} = \sum_{s \in V(G)^{|V(H)|}}ind(s\cong H)$. Similarly, $\binom{G}{H}$ refers to the number of subgraphs of $G$ isomorphic to $H$. \\
%Let $N'(v)$ refer to the subgraph induced by the neighborhood of $v$. 
%$K_{1,H}$ refers to $H + v$ with $v$ adjacent to all vertices in $H$. \\
$H \subseteq G$ refers to $H$ being a subgraph of $G$\\
$\smallSq{G}{H}$ represents the number of times $H$ appears as an induced subgraph of $G$.\\
$\binom{G}{H}$ represents the number of times $H$ appears as a subgraph of $G$.\\
$\mathcal{F}$ will always refer to a sequence of objects.

\section{Generalization}
 Kelly's lemma, a deeply important tool in reconstructibility, implies that $\smallSq{G}{H}$ and $\binom{G}{H}$ are reconstructible for any graph $H$ such that $|V(H)| < |V(G)|$ \cite{kelly1957}. The concept of a covering is used where we say $\mathcal{G} = (G_1, G_2,...G_k)$ is a covering of $H \subseteq G$ when and only when $\cup G_i = H$. $c(\mathcal{F},X)$ refers to the number of coverings of $X$ where $X$ is an unlabeled graph, and $\mathcal{F}$ is some sequence of unlabeled graphs such that, for each covering $\mathcal{G}$, $G_i \cong F_i$. This will be used in Kocay's Lemma, which states the following\cite{kocay1981}\cite{kocay1982}:
\[
\prod_\mathcal{F} \binom{G}{F_i} = \sum_{X \subseteq G} c(\mathcal{F},X) \binom{G}{X}
\]

The proof of Kocay's lemma follows nicely from a counting argument for both sides and serves as a powerful tool in graph reconstruction. Knowing the appearance of subgraphs with fewer vertices allows us to explicitly calculate information about only subgraphs with $n$ vertices. Beyond that knowledge even, it still provides information about several subgraphs using only information about a few. Kocay's lemma has been generalized in contexts of subgraph posets and such \cite{thatte2015}, but this paper explores a less specialized approach.

Kocay's lemma can be generalized from the following (where $f^{-1}(y)$ denotes the preimage of y).
\begin{lemma}
For any function $f: X\rightarrow Y$ and any $X, Y$: $|X| = \sum_{y\in Y} |f^{-1}(y)|$
\end{lemma}
Proof: Trivial $\square$\\

To apply this to Kocay's lemma, we will now begin to explain the notion of covering in a general sense.\\ Firstly, let $K$ be some set of objects such that each object has two components: a single component representing identity (call this the $L$ component), where $L$ is also the set of all $L$ components, and a secondary component representing multiplicity. For example, if I take the set of subgraphs of a graph with at least two squares, I can say there there exists $\{ (C_4,1),(C_4,2)\}$ within the set, where $C_4 \in L$ represents a square. For the multiplicity of an element in $H\in L$ within $K$, we will call it $\binom{K}{H}$. Say $G_i \cong G_j$ for $G_i, G_j \in K$ when and only when they have the same $L$ component, and we'll say $G_i \cong F_i$ for $F_i \in L$ when the $L$ component of $G_i$ is $F_i$.\\
Let's say $\mathcal{F}$ is a sequence of elements of $L$. Also, we'll say that a \textit{cover} of $B \in K$ by $\mathcal{F} = (F_1, F_2... F_k)$ is a sequence $(G_1, G_2...G_k)$ such that $G_{i}\cong F_i$, $G_i \in K$ and $\cup G_i = B$. To imagine this concretely, it is as though you have a bunch of substructures of $B$ whose union is exactly $B$ and such that each object in the sequence $(G_1, G_2...G_k)$ represents a component of $\mathcal{F}$. In this case, $c(\mathcal{F},B)$ is the number of covers of $B$ by $\mathcal{F}$. Say that $G$ is such that $\forall G_{i} \cong G_{j}: c(\mathcal{F},G_i) = c(\mathcal{F},G_j)$. This allows us to define $c(\mathcal{F}, X)$ for $X \in L$ as such: $c(\mathcal{F},X) = c(\mathcal{F},G_i): G_i \cong X \land G_i \in K$.
\begin{Theorem}
For any such system, $\prod_\mathcal{F}\binom{K}{F_i} =  \sum_{X\in L} c(\mathcal{F},X)\binom{K}{X}$
\end{Theorem}

Now, consider taking the set of all instances of each member of $\mathcal{F}$ and taking the cartesian product.\\
\[
\begin{pmatrix}
F_{1,1} \\
F_{1,2}\\
\vdots \\
F_{1,\binom{G}{F_1}}

\end{pmatrix}\times
\begin{pmatrix}
F_{2,1} \\
F_{2,2}\\
\vdots \\
F_{2,\binom{G}{F_2}}

\end{pmatrix}\times...
\times
\begin{pmatrix}
F_{|\mathcal{F}|,1} \\
F_{|\mathcal{F}|,2}\\
\vdots \\
F_{|\mathcal{F}|,\binom{G}{F_{|\mathcal{F}|}}}

\end{pmatrix}
\]
\\
Let's call this set $A$.
We can see, importantly, that the set of covers of $\mathcal{F}$ onto $g$ is the preimage of $g$ under the union operation on elements of $A$.
As such, we can use lemma 1 to acquire a simple equation around this.
\[
|A| = \sum_{g\in G} c(\mathcal{F},g)
\]
Next, it is clear that $|A| = \prod_\mathcal{F}\binom{G}{F_i}$, and so we can further refine the equation.\\
\[
\prod_\mathcal{F}\binom{G}{F_i} =  \sum_{g\in G} c(\mathcal{F},g)
\]
Next, we can apply the fact that $c(\mathcal{F},g_i) = c(\mathcal{F},g_j)$ for $g_i \cong g_j$, and partition this by isomorphism class (the set of which is $L$). We'll say $g \cong X \in L$.
\[
\prod_\mathcal{F}\binom{G}{F_i} =  \sum_{X\in L} c(\mathcal{F},X)\binom{G}{X}
\]
$\square$\\

This can be easily extended to functions with similar properties of congruence as coverings and unions like for discrete sets of points in space, but this will not be explored in this paper.\\

\section{Extension and Path Examples}
%We will use an association of finite simple graphs to a 2-edge colored complete graph. Suppose we have red and blue edges.
%So, we can define a labeled 2-edge colored graph $J$ by some $V(J)$, $R(J)$, $B(J)$ such that $R(J) \cap B(J) = \emptyset$ and $R(J) \cap B(J) = \binom{V(J)}{2}$.\\
We can now find some applications using Hamiltonian paths, known to be reconstructible but non-trivially \cite{tutte1979}. Firstly, the set of subgraphs under unions into connected graphs allows for Kocay's lemma, with all subgraphs. It does not pertain only to those with order $<n$.  Kocay's lemma is often used in conjunction with Kelly's lemma using $\mathcal{F}$ as some sequence of graphs $|V(F_i)|<n$, such as in \cite{kocay1982}, but it functionally can be used very powerfully with inputs of greater order. Below is an example.\\

Trivially, we can reconstruct $\binom{G}{\lceil \frac{n-1}{2}\rceil K_2}$ and $\binom{G}{\lfloor \frac{n-1}{2}\rfloor K_2}$ as they are disconnected for $n>4$. Unlike the traditional form of the lemma, we can assert that the lemma holds true for all $\mathcal{F}$ composed of subgraphs. They needn't have fewer vertices than $G$. As such, we can reconstruct the number of connected graphs of order n which they reconstruct, multiplied by their covering numbers. However, the only order n connected graph covered by $\lceil \frac{n-1}{2}\rceil K_2$ and $\lfloor \frac{n-1}{2}\rfloor K_2$ is $P_n$. This can be seen in that it would have n-1 edges and be connected. Now, suppose it had a vertex of degree 3. Since $\mathcal{F}$ has two components and has this vertex in its covers, it must cover this vertex. This implies that one member of $\mathcal{F}$ covers at least 2 of the edges. This would imply $P_2$ is a subgraph of one of the components of $\mathcal{F}$, which is not the case. As such, all graphs reconstructed are trees with $\Delta \leq 2$. This exactly describes a $P_n$. As such, $c(\mathcal{F}, P_n)\binom{G}{P_n}$ is a reconstructible quantity. Since $c(\mathcal{F})$ is calculable, $\binom{G}{P_n}$ is reconstructible. $\square$\\

Next, we will  rigorously develop the tools for an extension of Kocay's lemma into 2-edge refined graphs, such as the one used in his work to reconstruct $\binom{G}{P_n}$\cite{kocay1982}.

For any finite, simple graph $G$, we can associate a canonical form (ignoring color) of 2-edge refined complete graph by associating $G$ with $G'$ as follows:\\
\[
V(G') = V(G)
\]
\[
R(G') = E(G)
\]
\[
B(G') = \binom{V(G)}{2} - R(G')
\]
This can be seen as just taking all non-edges as blue edges, and all edges as red edges. Now, we need an important step.\\
\begin{lemma} [Weinstein. 1975 \cite{weinstein1975}]
$\smallSq{G'}{H}$ and $\binom{G'}{H}$ are reconstructible for $H$ with $|V(H)| < |V(G')|$
\end{lemma}

This is very useful already, but it also allows us to reconstruct interesting subgraphs. The subgraphs which are reconstructible as in Kelly's lemma are simply 2-edge refined graphs, but with an interesting caveat.\\
Say $H \subseteq G'$.
Then $V(H) \subseteq V(G')$, $R(H) \subseteq R(G')$, $B(H) \subseteq B(G')$ and $B(H), R(H) \subseteq V(H) \times V(H)$. 
If we try to represent this in G, it is similar to a subgraph, but it preserves adjacency and nonadjacency. This allows us to do some nice and funky things. Firstly, let $\mathcal{F}$ represent any sequence of isomorphism classes of subgraphs of $G'$.
\begin{Theorem}
$\prod_\mathcal{F}\binom{G'}{F_i} =  \sum_{X\in L} c(\mathcal{F},X)\binom{G'}{X} =\sum_{X \subseteq G'} c(\mathcal{F},X)\binom{G'}{X}$
\end{Theorem}
Proof: We can establish this with theorem 1. Firstly, we can take the set of isomorphism classes of subgraphs as $L$ (as in theorem 1), and the set of subgraphs with multiplicity as $K$. For 2-edge refined graphs, we can see that $c(\mathcal{F}, G'_i) = c(\mathcal{F},G'_j)$ when $G'_i \cong G'_j$. There must be a bijection, say $\varphi:V(G'_i)\rightarrow V(G'_j)$ which preserves both red and blue adjacency as well as nonadjacency by the definition of congruence on 2-edge refined graphs. That is,
If \( G'_i \cong G'_j \), then there exists a bijection \( \varphi : V(G'_i) \to V(G'_j) \) such that for all \( u,v \in V(G'_i) \), the following hold:
\begin{itemize}
    \item \( \{u, v\} \in R(G'_i) \iff \{\varphi(u), \varphi(v)\} \in R(G'_j) \),
    \item \( \{u, v\} \in B(G'_i) \iff \{\varphi(u), \varphi(v)\} \in B(G'_j) \),
    \item \( \{u, v\} \notin R(G'_j) \cup B(G'_j) \iff \{\varphi(u), \varphi(v)\} \notin R(G'_j)\cup B(G'_j) \)
\end{itemize}
Now, say we take any covering $\mathcal{G}_1$ of $G'_i$ (regardless of $\mathcal{F}$). This corresponds to a sequence of 2-edge refined graphs, $(g_1, g_2...g_{|\mathcal{F}|})$. Each of these is exactly definable as a subset of the colored edges and vertices of $G_i$. Say we take any element in the covering sequence, let's call it $g_k$. As each edge in $g_k$ is in $G_i$, as well as every vertex and nonedge, we have:
\begin{itemize}
   \item \( \{u, v\} \in R(g_k) \iff \{\varphi(u), \varphi(v)\} \in R(G'_j) \),
    \item \( \{u, v\} \in B(g_k) \iff \{\varphi(u), \varphi(v)\} \in B(G'_j) \),
\end{itemize}
Knowing such, say we take the subgraph $g'_k$ of $G'_j$ with $V(g'_k) = \varphi(V(g_k))$ and all the edges from $\varphi(g_k)$. We then have a $\varphi: V(g_k) \rightarrow V(g'_k)$ such that all colored adjacencies and nonadjacencies are preserved, giving us a bijection. As such, we can associate a bijection between the elements of each cover, meaning the number of covers is preserved across the isomorphism class of a 2-edge refined graph. $\square$\\

Now, we can do a similar proof using an edge coloring of $K_n$. Firstly, we reconstruct the canonical 2-coloring of the deck of $G$, and then we can proceed.\\

We can apply theorem 2 and the same reasoning that the number of 2-connected graphs with fixed number of edges are reconstructible. \\
\begin{lemma}
$\sum_{X\subseteq G': |V(X)|=n} c(\mathcal{F},X)\binom{G'}{X}$ is reconstructible for $|V(F_i)|<n$
\end{lemma}
Proof: 
We know $D(G')$ and all induced subgraphs of $G'$. As such, we can reconstruct $\prod_{\mathcal{F}} \binom{G'}{F_i}$ for $|V(F_i)|<n$. By theorem 2, we can reconstruct  \\
$\sum_{X\subseteq G'} c(\mathcal{F},X)\binom{G'}{X}$.\\
However, analogously to Kocay's use for subgraphs, we can calculate \\
$\sum_{X\subseteq G': |V(X)|< n} c(\mathcal{F},X)\binom{G'}{X}$ using Kelly's lemma for $\binom{G'}{X}$ and just calculating $c(\mathcal{F}, X)$. Then, we can subtract such from the original sum to get \\
$\sum_{X\subseteq G': |V(X)|=n} c(\mathcal{F},X)\binom{G'}{X}$. $\square$

Now, we show that we can also reconstruct $\sum_{X\in D} c(\mathcal{F},X)\binom{G'}{X}$ where $D$ is the set of spanning disconnected subgraphs for $\mathcal{F}$ reconstructible. 

\begin{lemma}
$\sum_{X\in D} c(\mathcal{F},X)\binom{G'}{X}$ where $D$ is the set of isomorphism classes of spanning disconnected subgraphs for $\mathcal{F}$ reconstructible.
\end{lemma}

Proof: By lemma 3, we can reconstruct $\sum_{X\subseteq G': |V(X)|=n} c(\mathcal{F},X)\binom{G'}{X}$. Now, say we wish to reconstruct $\sum_{X\in D} c(\mathcal{F},X)\binom{G'}{X}$ where $D$ is the set of isomorphism classes of spanning disconnected subgraphs of $G'$. Let $D'$ be the set of all disconnected 2-edge refined graphs on $n$ vertices. $\binom{G'}{d'}$ for $d' \in D'$ is reconstructible, as we can simply reconstruct $\sum_{X\subseteq G': |V(X)|=n} c(\mathcal{G},X)\binom{G'}{X}$where $\mathcal{G}$ %
is  a sequence of the connected components of $d'$, and then we can simply divide by $c(\mathcal{G},d')$ to get $\binom{G'}{d'}$. Next, we can calculate $c(\mathcal{F},d')$ for all $d' \in D'$ and build the sum $\sum_{d'\in D'} c(\mathcal{F},d')\binom{G'}{d'}$. $c(\mathcal{F},d')$ will be 0 when $d'$ is not a spanning subgraph of $G'$, so we can rewrite this sum as $\sum_{X\in D} c(\mathcal{F},X)\binom{G'}{X}$. $\square$\\
Also notice how it follows from the traditional proof for graphs and subgraphs. \\

The same process follows for 1-connected 2-edge refined subgraphs which have blocks $\mathcal{G}$. That is, 
$\sum_\mathcal{G} |V(G_i)| = n-|\mathcal{G}|+1$ and any two components can be separated by the removal of a single vertex., as we can count the compositions of these and remove any disconnected 2-edge refined subgraphs\\

It is true from such that we can reconstruct the number of connected 2-edge refined subgraphs with $n-1$ edges (number of red + blue), as we can simply use lemma 3 and subtract all the disconnected results.\\
As such, we can use the modified Kocay's lemma to reconstruct the number of 2-connected graphs with $n-1$ red edges and a single blue edge, by letting such be $\mathcal{F}$ and removing all the 1-connected and disconnected results.\\
Then, supposing we do this up to $k$ edges, we can remove all the disconnected, 1-connected graphs, and 2-connected graphs with fewer than k edges. $\square$

Now, we have shown that each tool can be easily redeveloped from Kocay's original lemma to prove his result: Every Hamiltonian path is in exactly 1 red Hamiltonian cycle or 1 Hamiltonian cycle with only a single blue edge. Since these are both reconstructible as well as the number of cycles in each, $\binom{G}{P_n}$ is reconstructible. \\

This proof corresponds very closely to the Tutte's original, except instead of indirectly counting the lack of an edge, this directly counts it: circumventing the need for a large number of linear equations and multiedges entirely. Notably, that solution also implicitly uses a modified form of Kocay's lemma.\\
While it may seem that it took a bit of work, the tools developed are standard for the typical application of Kocay's lemma.\\

The reconstruction conjecture is equivalent to supposing that a graph can be reconstructed using only sums from Kocay's lemma.

As such, no result presented on subgraph counts is unattainable theoretically by (potentially repeated) application of Kocay's lemma. However, the repeated use is not known to necessarily be restricted to finite repetition for general graphs, and the finite results may also be significantly harder to attain regardless.\\

Beyond that, though, there are also novel results that are easily attained with a modified form of Kocay's lemma.\\

\section{Result on Trees}

For any tree T and any graph G, we have the following:
\begin{Theorem}
 $\binom{G}{T}$, $\binom{G}{T} + \binom{\overline{G}}{T}$ or $\binom{G}{T} - \binom{\overline{G}}{T}$ is reconstructible.
\end{Theorem}

Proof: Firstly, we use $H$ as some 2-edge colored $K_n$. We will denote an edge $e_i$ and specify the color like $be_i$ for blue, and $re_i$ for red. $H_{be_i}$ will represent $H$ but with $e_i$ as blue, and the same for red. We will take the 2-form of $G$.

For the following lemma, consider $e_i$ to be $(v_i,u_i)$
\begin{lemma} 
$\binom{G}{H-e_i}*|orb_{H-e_i}(v_i,u_i)|= \binom{G}{H_{be_i}}*\binom{H_{be_i}}{H-e_i} + \binom{G}{H_{re_i}}*\binom{H_{re_i}}{H-e_i}$
\end{lemma}
Proof: We are going to try to count $H_{be_i}$ and $H_{re_i}$ to illustrate. Consider every instance of $H-e_i$ in $G$. As $G$ is a complete graph, there must be some edge present in place of $e_i$. So, each of these must correspond to an $H_{be_i}$ or an $H_{re_i}$. The total number counted by each instance of $H-e_i$ is $|orb_{H-e_i}(v_i,u_i)|$. For every instance of $H_{be_i}$ counted, however, it is counted exactly $\binom{H_{be_i}}{H-e_i}$ times. The same applies for $H_{re_i}$, and so we get the above equation. $\square$.\\

Now, we can modify this for our advantage. Say we select any specific edge $e=(v,u)$ in some 2-edge coloring of $T$. Then, we can say that $T-(v,u)$ is reconstructible, as it is disconnected. Since $|orb_{T-e}(v_i,u_i)|$ is calculable, by the previous lemma, we can reconstruct $ \binom{G}{T_{be}}*\binom{T_{be}}{T-e} + \binom{G}{T_{re}}*\binom{T_{re}}{T-e}$. Now, we can proceed with a sort of descent argument. \\
For the sake of convenience, order the edges in an entirely red $T$. Let's say the sequence is $S=(e_1, e_2, e_3... e_{n-1})$. 
Again, for the sake of brevity, let's say the nth position is $S_n$ and $T_{bS_n}$ refers to $T$ with all edges red except for blue edges up to $S_n$.\\
We can then reconstruct $\binom{G}{T_{bS_i}}*\binom{T_{bS_i}}{T_{bS_i}-e_{i+1}}+\binom{G}{T_{bS_{i+1}}}*\binom{T_{bS_{i+1}}}{T_{bS_{i+1}}-e_{i+1}}$. Notice that each instance of such is a linear combination of $\binom{G}{T_{bS_i}}$ and $\binom{G}{T_{bS_{i+1}}}$. As such, using such linear combinations from $S_1$ to $S_{n-1}$, we can acquire some linear combination of $\binom{G}{T}$ and $\binom{G}{T_{bS_{n-1}}}$. Note: $\binom{G}{T_{bS_{n-1}}} = \binom{\overline{G}}{T}$.

This is already somewhat nice, but consider if that linear combination is not $\binom{G}{T} + \binom{\overline{G}}{T}$ or $\binom{G}{T} - \binom{\overline{G}}{T}$ or some scalar multiple of either. Notice, all the coefficients of these terms are independent of $G$. But beyond that, we can notably reconstruct the deck of the complement of the graph. As such, we can repeat the process with $\overline{G}$ to achieve the same equation with the coefficients reversed. That is, if we can reconstruct  $a\binom{G}{T} + b\binom{\overline{G}}{T}$, then $a\binom{\overline{G}}{T} + b\binom{\overline{\overline{G}}}{T} = b\binom{G}{T} + a\binom{\overline{G}}{T}$ is reconstructible. As such, either $a = \pm b$, or these equations can derive $\binom{G}{T}$. $\square$
\section{Conclusion}
We have rigorously developed the 2-edge refined form of Kocay's lemma, as well as showed that it can be generalized to many different spaces beyond graphs, as well as showing a novel result for the reconstruction of the number of trees in a graph.
\bibliographystyle{plain} % or alpha, ieeetr, apalike, etc.
\bibliography{references.bib}

\begin{thebibliography}{1}

\bibitem{kelly1957}
P.~J. Kelly.
\newblock A congruence theorem for trees.
\newblock {\em Pacific Journal of Mathematics}, 7:961--968, 1957.

\bibitem{kocay1981}
W.~L. Kocay.
\newblock On reconstructing spanning subgraphs.
\newblock {\em Ars Combinatoria}, 11:301--313, 1981.

\bibitem{kocay1982}
W.~L. Kocay.
\newblock Some new methods in reconstruction theory.
\newblock In E.~J. Billington, S.~Oates-Williams, and A.~Penfold Street,
  editors, {\em Combinatorial Mathematics IX}, volume 952 of {\em Lecture Notes
  in Mathematics}, pages 89--114. Springer-Verlag, Berlin, 1982.
\newblock (Proc. 9th Australian Conf. on Combinatorial Mathematics, Univ. of
  Queensland, Brisbane).

\bibitem{thatte2015}
Bhalchandra~D. Thatte.
\newblock Subgraph posets and graph reconstruction, 2015.

\bibitem{tutte1979}
W.~T. Tutte.
\newblock All the king’s horses—a guide to reconstruction.
\newblock In J.~A. Bondy and U.~S.~R. Murty, editors, {\em Graph Theory and
  Related Topics}. Academic Press, 1979.

\bibitem{weinstein1975}
Joseph~M. Weinstein.
\newblock Reconstructing colored graphs.
\newblock {\em Pacific Journal of Mathematics}, 57(1):307--314, 1975.

\end{thebibliography}

\end{document}